\input amstex
\documentstyle{amsppt}
\magnification=\magstep1                        
\hsize6.5truein\vsize8.9truein                  
\NoRunningHeads
	\loadeusm

\magnification=\magstep1                        
\hsize6.5truein\vsize8.9truein                  
\NoRunningHeads
\loadeusm

\document
\topmatter

\title
Improved lower bounds for the Mahler measure of the Fekete polynomials
\endtitle

\rightheadtext{the Mahler measure of the Fekete polynomials}

\author Tam\'as Erd\'elyi
\endauthor

\address Department of Mathematics, Texas A\&M University,
College Station, Texas 77843, College Station, Texas 77843 (T. Erd\'elyi) \endaddress

\endthanks


\date February 17, 2017
\enddate

\abstract
We show that there is an absolute constant $c > 1/2$ such that the Mahler measure of the Fekete polynomials 
$f_p$ of the form
$$f_p(z) := \sum_{k=1}^{p-1}{\left( \frac kp \right)z^k}\,,$$
(where the coefficients are the usual Legendre symbols) 
is at least $c\sqrt{p}$ for all sufficiently large primes $p$. 
This improves the lower bound $\left(\frac 12 - \varepsilon\right)\sqrt{p}$ known before for 
the Mahler measure of the Fekete polynomials $f_p$ for all sufficiently large primes 
$p \geq c_{\varepsilon}$. Our approach is based on the study of the zeros of the 
Fekete polynomials on the unit circle.
\endabstract

\endtopmatter

\head 1. Introduction and Notation \endhead

Let $D$ be the open unit disk of the complex plane. Its boundary, the unit circle
of the complex plane, is denoted by $\partial D$. Let
$${\Cal K}_n := \left\{P_n: P_n(z) = \sum_{k=0}^n{a_k z^k}, \enskip  a_k \in {\Bbb C}\,,\enskip |a_k| = 1 \right\}\,.$$
The class ${\Cal K}_n$ is often called the collection of all (complex) unimodular
polynomials of degree $n$. Let
$${\Cal L}_n := \left\{P_n: P_n(z) = \sum_{k=0}^n{a_k z^k}, \enskip  a_k \in \{-1,1\} \right\}\,.$$ 
The class ${\Cal L}_n$ is often called the collection of all Littlewood
polynomials of degree $n$. By Parseval's formula,
$$\int_{0}^{2\pi}{|P_n(e^{it})|^2 \, dt} = 2\pi(n+1)$$
for all $P_n \in {\Cal K}_n$. Therefore
$$\min_{z \in \partial D}{|P_n(z)|} < \sqrt{n+1} < \max_{z \in \partial D}{|P_n(z)|}$$
for all $P_n \in {\Cal K}_n$ and $n \geq 1$.
An old problem (or rather an old theme) is the following.

Let $\alpha < \beta$ be real numbers. The Mahler measure $M_{0}(P,[\alpha,\beta])$ is
defined for bounded measurable functions $P(e^{it})$ defined on $[\alpha,\beta]$ as
$$M_{0}(P,[\alpha,\beta]) := \exp\left(\frac{1}{\beta - \alpha} \int_{\alpha}^{\beta}
{\log|P(e^{it})|\,dt} \right)\,.$$
It is well known that
$$M_{0}(P,[\alpha,\beta]) = \lim_{q \rightarrow 0+}{M_{q}(P,[\alpha,\beta])}\,,$$
where, for $q > 0$,
$$M_{q}(P,[\alpha,\beta]) := \left( \frac{1}{\beta-\alpha} \int_{\alpha}^{\beta}
{\left| P(e^{it}) \right|^q\,dt} \right)^{1/q}\,.$$
It is a simple consequence of the Jensen formula that
$$M_0(P) := M_0(P,[0,2\pi]) = |c| \prod_{k=1}^n{\max\{1,|z_k|\}}$$
for every polynomial of the form
$$P(z) = c\prod_{k=1}^n{(z-z_k)}\,, \qquad c,z_k \in {\Bbb C}\,.$$
P. Borwein and Lockhart [B-01]  investigated the asymptotic behavior of the mean
value of normalized $L_q$ norms of Littlewood polynomials for arbitrary $q > 0$.
Using the Lindeberg Central Limit Theorem and dominated convergence, they proved that
$$\lim_{n \rightarrow \infty} {\frac{1}{2^{n+1}} \sum_{f \in {\Cal L}_n}
{\frac{(M_q(f,[0,2\pi]))^q}{n^{q/2}}}}= \Gamma \left( 1+ \frac q2 \right)$$
for every $q > 0$. In [C-15] we proved that
$$\lim_{n \rightarrow \infty} {\frac{1}{2^{n+1}} \sum_{f \in {\Cal L}_n}
{\frac{M_q(f,[0,2\pi])}{n^{1/2}}}} = \left( \Gamma \left( 1+ \frac q2 \right) \right)^{1/q}$$
for every $q > 0$. We also proved analogous results for the Mahler measure. Namely, 
using the notation $\widehat{f}(z) := \max\{|f(z)|,n^{-1}\}$, we have
$$\lim_{n \rightarrow \infty} {\frac{1}{2^{n+1}} \sum_{f \in {\Cal L}_n}
{\log \left( \frac{M_0(\widehat{f},[0,2\pi])}{n^{1/2}}\right)}} = -\gamma/2\,,$$ 
and
$$\lim_{n \rightarrow \infty} {\frac{1}{2^{n+1}} \sum_{f \in {\Cal L}_n}
{\frac{M_0(f,[0,2\pi])}{n^{1/2}}}} = e^{-\gamma/2}\,,$$
where
$$\gamma := \lim_{n \rightarrow \infty}{\left( \sum_{k=1}^n{\frac 1k - \log n} \right)} = 0.577215 \ldots$$
is the Euler constant and $e^{-\gamma/2} = 0.749306\ldots$.
These are analogues of the results proved earlier by Choi and Mossinghoff
[C-11] for polynomials in ${\Cal K}_n$.

Finding polynomials with suitably restricted coefficients and maximal Mahler measure 
has interested many authors. Beller and Newman [B-73] constructed unimodular polynomials of 
degree $n$ whose Mahler measure is at least $\sqrt{n}-c/\log n$.
For a prime $p$ the $p$-th Fekete polynomial is defined as
$$f_p(z) := \sum_{k=1}^{p-1}{\left( \frac kp \right)z^k}\,,$$
where
$$\left( \frac kp \right) =
\cases
1, \quad \text{if \enskip} x^2 \equiv k \enskip (\text {mod\,}p) \enskip
\text{for an} \enskip x \neq 0\,,
\\
0, \quad \text{if \enskip} p \enskip \text{divides} \enskip k\,,
\\
-1, \quad \text{otherwise}
\endcases$$
is the usual Legendre symbol. Note that $g_p(z) := f_p(z)/z$ is a Littlewood
polynomial, and has the same Mahler measure as $f_p$.

Montgomery [M-80] proved the following fundamental result.

\proclaim{Theorem 1.1}
There are absolute constants $c_1 > 0$ and $c_2 > 0$ such that
$$c_1 \sqrt{p} \log \log p \leq \max_{z \in \partial D}{|f_p(z)|} \leq c_2 \sqrt{p} \log p\,.$$
\endproclaim

In [E-07] we proved the following result.

\proclaim{Theorem 1.2} For every $\varepsilon > 0$ there is
a constant $c_{\varepsilon}$ such that
$$M_0(f_p,[0,2\pi]) \geq \left(\frac 12 - \varepsilon \right)\sqrt{p}$$
for all primes $p \geq c_{\varepsilon}$.
\endproclaim

From Jensen's inequality,
$$M_{0}(f_{p}) \leq M_2(f_p) = \sqrt{p-1}\,.$$
However, as it was observed in [E-07], $\frac{1}{2}-\varepsilon$ in Theorem 1.2 cannot be replaced by 
$1-\varepsilon$. Indeed, if $p \geq 3$ is a prime and $p=4m+1$, then 
$f_p$ is self-reciprocal, that is, $z^pf_{p}(1/z) = f_{p}(z)$,
and hence
$$f_{p}(e^{2it}) = e^{ipt}\sum_{k=0}^{(p-3)/2}a_k \cos((2k+1)t) ,\qquad a_{k} \in \{-2,2\}\,.$$ 
Therefore a result of Littlewood [L-66] implies that
$$\split M_{0}(f_p) \leq & \, \frac{1}{2\pi} \int_{0}^{2\pi}{|f_{p}(e^{iu})| \, du} = 
\frac{1}{2\pi} \int_{0}^{\pi}{|f_{p}(e^{2it})| \, 2dt} = \frac{1}{2\pi} \int_{0}^{2\pi}{|f_{p}(e^{2it})| \, dt} \cr 
\leq & \, (1-\varepsilon_0) \sqrt{p-1} \cr \endsplit$$
with some absolute constant $\varepsilon_0 > 0$. 
If $p \geq 3$ is a prime and $p=4m+3$, then $f_{p}$ is anti-self-reciprocal, that is,
$z^pf_p(1/z) = -f_{p}(z)$,
and hence
$$if_{p}(e^{2it}) =e^{ipt}\sum_{k=0}^{(p-3)/2}a_k \sin((2k+1)t) ,\qquad a_{k} \in \{-2,2\}\,.$$ 
Therefore a result of Littlewood [L-66] implies that
$$\split M_{0}(f_p) & \, \leq \frac{1}{2\pi} \int_{0}^{2\pi}{|if_{p}(e^{iu})| \, du} = 
\frac{1}{2\pi} \int_{0}^{\pi}{|if_{p}(e^{2it})| \, 2dt} = \frac{1}{2\pi} \int_{0}^{2\pi}{|if_{p}(e^{2it})| \, dt} \cr  
& \, \leq (1-\varepsilon_0) \sqrt{p-1} \cr \endsplit$$
with some absolute constant $\varepsilon_0 > 0$.

It is an interesting open question whether there is a sequence of Littlewood polynomials $(f_n)$ such that for an arbitrary 
$\varepsilon >0$, and $n$ large enough,
$$M_{0}(f_n) \geq (1-\varepsilon) \sqrt{n}\,.$$ 

In [E-11] Theorem 1.2 was extended to subarcs of the unit circle.

\proclaim{Theorem 1.3} There exists an absolute constant $c_1 > 0$ such that
$$M_0(f_p,[\alpha,\beta]) \geq c_1 p^{1/2}$$
for all primes $p$ and for all $\alpha, \beta \in {\Bbb R}$ such that
$(\log p)^{3/2}p^{-1/2} \leq \beta - \alpha \leq 2\pi$. 
\endproclaim

In [E-12] we gave an upper bound for the average value of $|f_p(z)|^q$ over any
subarc $I$ of the unit circle, valid for all sufficiently large primes $p$ and all
exponents $q > 0$.

\proclaim{Theorem 1.4} There exists a constant $c_2(q,\varepsilon)$ depending only on $q > 0$
and $\varepsilon > 0$ such that
$$M_q(f_p,[\alpha,\beta]) \leq c_2(q,\varepsilon)p^{1/2}\,,$$
for all primes $p$  and for all $\alpha, \beta \in {\Bbb R}$ such that
$2p^{-1/2+\varepsilon} \leq \beta - \alpha \leq 2\pi$.
\endproclaim

We remark that a combination of Theorems 1.3 and 1.4 shows that there is an absolute
constant $c_1 > 0$ and a constant $c_2(q,\varepsilon) > 0$ depending only on $q > 0$
and $\varepsilon > 0$ such that
$$c_1p^{1/2} \leq M_q(f_p,[\alpha,\beta])  \leq c_2(q,\varepsilon)p^{1/2}$$
for all primes $p$  and for all $\alpha, \beta \in {\Bbb R}$ such that
$(\log p)^{3/2}p^{-1/2} \leq 2p^{-1/2+\varepsilon} \leq \beta - \alpha \leq 2\pi$.

The $L_q$ norm of polynomials related to Fekete polynomials were studied in several 
recent papers. See [B-01b], [B-02], [B-04], [G-16], [J-13a], and [J-13b], for example. 
An interesting extremal property of the Fekete polynomials is proved in [B-01c]. 

Fekete might have been the first one to study analytic properties of the Fekete polynomials. 
He had an idea of proving non-existence of Siegel zeros (that is, real zeros 
``especially close to $1$") of Dirichlet $L$-functions 
from the positivity of Fekete polynomials on the interval $(0,1)$, 
where the positivity of Fekete polynomials is often referred to as 
the Fekete Hypothesis. There were many mathematicians 
trying to understand the zeros of Fekete polynomials including 
Fekete and P\'olya [F-12], P\'olya [P-19], Chowla [C-35], Heilbronn [H-37], 
Montgomery [M-80], Baker and Montgomery [B-90], and Jung and Shen [J-16].  

Baker and Montgomery [B-90] proved that $f_p$ has a large number of zeros 
in $(0,1)$ for almost all primes $p$, that is, the number of zeros of 
$f_p$ in $(0,1)$ tends to $\infty$ as $p$ tends to $\infty$, and it 
seems likely that there are, in fact, about $\log\log p$ such zeros.  

Conrey, Granville, Poonen, and Soundararajan [C-00] showed that $f_p$ has
asymptotically $\kappa p$ zeros on the unit circle, where $0.500668 < \kappa < 0.500813$.

An interesting recent paper [B-17] studies power series approximations to Fekete polynomials. 

It is conjectured, see [B-02] for instance, that there are sequences of
flat Littlewood polynomials $P_n \in {\Cal L}_n$ satisfying
$$c_1 \sqrt{n+1}  \leq |P_n(z)| \leq c_2 \sqrt{n+1}\,, \qquad z \in \partial{D}\,,$$
with absolute constants $c_1>0$ and $c_2>0$. However, the lower bound part 
of this conjecture, by itself, seems hard, and no sequence is known that satisfies 
just the lower bound. A sequence of Littlewood polynomials satisfying just the upper 
bound is given by the Rudin-Shapiro polynomials. They appear in Harold Shapiro's 1951 
thesis [S-51] at MIT and are sometimes called just Shapiro polynomials. They also arise 
independently in a paper by Golay (1951). They are remarkably simple to construct and 
are a rich source of counterexamples to possible conjectures. The Rudin-Shapiro 
polynomials are defined recursively as follows:
$$\split P_0(z) & :=1\,, \qquad Q_0(z) := 1\,, \cr 
P_{n+1}(z) & := P_n(z) + z^{2^n}Q_n(z)\,, \cr
Q_{n+1}(z) & := P_n(z) - z^{2^n}Q_n(z)\,, \qquad n = 0,1,2,\ldots \,. \cr \endsplit$$
Note that both $P_n$ and $Q_n$ are polynomials of degree $N-1$ with $N := 2^n$ having each
of their coefficients in $\{-1,1\}$.
In [E-16] we showed that the Mahler measure and the maximum norm of the
Rudin-Shapiro polynomials on the unit circle of the complex plane have the
same size.

\proclaim{Theorem 1.5}
Let $P_n$ and $Q_n$ be the $n$-th Rudin-Shapiro polynomials defined in Section 1.
There is an absolute constant $c_1 > 0$ such that
$$M_0(P_n,[0,2\pi]) = M_0(Q_n,[0,2\pi]) \geq c_1\sqrt{N}\,,$$
where
$$N := 2^n = \text{\rm deg}(P_n)+1  = \text{\rm deg}(Q_n)+1\,.$$
\endproclaim

\head 2. New Result \endhead

In this paper we improve the factor $\left(\frac 12 - \varepsilon \right)$ in Theorem 1.1 to an absolute 
constant $c > 1/2$. Namely we prove the following. 

\proclaim{Theorem 2.1} There is an absolute constant $c > 1/2$ such that
$$M_0(f_p) \geq c\sqrt{p}$$
for all sufficiently large primes.
\endproclaim

\head 3. Lemmas \endhead

To prove the theorem we need a few lemmas. For a natural number $p$ let 
$$\zeta_p := \exp\left(\frac{2\pi i}{p}\right)$$
be the first $p$-th root of unity. Our first lemma
formulates an characteristic property of the Fekete polynomials.
A simple proof is given in [B-02, pp. 37-38].

\proclaim{Lemma 3.1 (Gauss)} We have
$$f_p(\zeta_p^j) = \sqrt{\left(\frac{-1}{p}\right)p}\,, \qquad j=1,2,\ldots,p-1\,,$$
and $f_p(1) = 0$.
\endproclaim

\proclaim{Lemma 3.2} We have
$$\left(\prod_{j=0}^{p-1}{|Q(\zeta_p^j)|}\right)^{1/p} \leq 2M_0(Q)$$
for all polynomials $Q$ of degree at most $p$ with complex coefficients.
\endproclaim

\proclaim{Lemma 3.3} Let $0 < \eta \leq \pi/2$ be fixed. Suppose
a polynomial $Q$ of degree at most $p$ with complex coefficients has at least $k$ zeros
$$b_j = e^{it_j}, \qquad j=1,2,\ldots,k\,,$$
such that
$$t_j \in [0,2\pi) \setminus \bigcup\limits_{\nu=0}^{p-1}
{\left(\frac{(2\nu+1)\pi}{p} - \frac{\eta}{p}, \enskip \frac{(2\nu+1)\pi}{p} + \frac{\eta}{p}\right)}\,.$$
We have
$$\left(\prod_{j=0}^{p-1}{|Q(\zeta_p^j)|}\right)^{1/p} \leq 2\left(\cos \frac{\eta}{2} \right)^{k/p}M_0(Q)\,.$$
\endproclaim

In the proof of Theorem 2.1 we need one of the following two results. For proofs see 
[B-97a] and [B97-b], respectively.

\proclaim{Lemma 3.4} There is an absolute constant $c > 0$ such that every
$Q \in {\Cal K}_n$ has at most $c\sqrt n$ real zeros.
\endproclaim

\proclaim{Lemma 3.5} There is an absolute constant $c > 0$ such that every
$Q \in {\Cal L}_n$ has at most $\displaystyle{\frac{c\log^2 n}{\log\log n}}$ zeros at $1$.
\endproclaim

The large sieve of number theory [M-78] asserts the following. 

\proclaim{Lemma 3.6} 
If
$$P(z) = \sum_{k=-n}^{n}{a_{k}z^{k}}\,, \qquad a_k \in {\Bbb C}\,,$$
is a trigonometric polynomial of degree at most $n$,
$$0 \leq t_1 < t_2 < \cdots < t_m \leq 2\pi\,,$$
and
$$\delta :=\min \left\{t_2 - t_1, t_3 - t_2,\dots,
t_m - t_{m-1}, 2\pi -\left(t_m - t_1 \right) \right\}\,,$$
then
$$\sum_{j=1}^{m}\left| P\left( e^{it_j} \right) \right|^{2}
\leq \left(\frac{2n+1}{2\pi}+\delta ^{-1}\right)
\int_{0}^{2\pi }{\left| P\left( e^{it}\right) \right|^{2}\,dt}\,.$$
\endproclaim

It turns out to be fairly easy to show that at least half of the zeros of $f_p$ are
on the unit circle $\partial D$. First note that
$$F_p(z) := z^{-p/2}f_p(z) = \sum_{a=1}^{(p-1)/2}{\left( \frac ap \right) \left( z^{a-p/2} + \left( \frac{-1}{p} \right) z^{p/2-a} \right)}\,.$$
Observe also that
$$F_p(e^{2i\pi t}) := \cases 
2 \displaystyle{\sum_{a=1}^{(p-1)/2}{\left( \frac ap \right)} \cos((2a-p)\pi t)} \quad \text{if} \enskip p \equiv 1 \enskip \text{mod\,} 4\\ 
2i \displaystyle{\sum_{a=1}^{(p-1)/2}{\left( \frac ap \right)} \sin((2a-p)\pi t)} \quad \text{if} \enskip p \equiv 3 \enskip \text{mod\,} 4\,.\\
\endcases \tag 3.1$$
Define $H_p(t) := F_p(e^{2i\pi t})$ if $p \equiv 1$ mod $\,4$, and $H_p(t) := -iF_p(e^{2i\pi t})$ if $p \equiv 3$ mod $\,4$.
By (3.1) we see that $H_p(t)$ is a periodic, continuous, real-valued function when $t$ is real.

\proclaim{Lemma 3.7} Let $p$ be a prime. There are at least $(p-3)/2$ values of $k \in \{0,1,\ldots,p-1\}$ for which
$H_p$ has a zero between $k/p$ and $(k+1)/p$.
\endproclaim

Our next lemma is Theorem 4 in [C-00]. For a proof of Lemma 3.8 below see Section 6 in [C-00].

\proclaim{Lemma 3.8}
Let $p$ be a prime. For every fixed real number $\delta$
$$\left| \left\{k \in \{1,2,\ldots,p\}: H_p \left( \frac{k+1/2}{p} \right) < \delta \sqrt{p} \right\} \right| \sim c_{\delta}p$$
as $p \rightarrow \infty$, where
$$c_{\delta} = \frac 12 + \frac{1}{\pi} \int_0^{\infty} {\sin(\delta \pi x) C(x) \, \frac{dx}{x}}\,, \qquad 
C(x) := \prod_{k=0}^{\infty}{\cos^2\left(\frac{2x}{2k+1}\right)}\,.$$
Moreover $c_{-\delta} = 1-c_{\delta}$ for all $\delta > 0$.
\endproclaim

\proclaim{Lemma 3.9}
For every $\varepsilon > 0$ there is a $\delta > 0$ such that
$$\left| \left\{k \in \{1,2,\ldots,p\}: H_p \left( \frac{k+1/2}{p} \right) \geq \delta \sqrt{p} \right\} \right| \geq (1-\varepsilon)p$$
for all sufficiently large primes $p \geq N_{\varepsilon}$.
\endproclaim

\proclaim{Lemma 3.10} Let $\gamma > 0$ be a real number. Let the subarcs $I_k$ of the unit circle 
$\partial D$ be defined by
$$I_k := \left\{e^{it}: \left|t - \frac{(2k+1)\pi}{p}\right| \leq \frac{\pi}{2p}\right\}\,, \qquad k=0,1,\ldots,p-1\,.$$
We have
$$m := |\{k \in \{0,1,\ldots,p-1\}: \max_{z \in I_k} {|f_p^{\prime}(z)|} \geq \gamma p^{3/2}\}| 
\leq \gamma^{-2}p/2$$
for all primes $p \geq 3$.
\endproclaim

\proclaim{Lemma 3.11} Given $\eta > 0$ let the subarcs $I_{k,\eta}$ of the unit circle $\partial D$ be defined by
$$I_{k,\eta} := \left\{e^{it}: \left|t - \frac{(2k+1)\pi}{p}\right| < \frac{\eta}{p}\right\}\,, \qquad k=0,1,\ldots,p-1\,.$$
For every $\varepsilon > 0$ there is an $\eta > 0$ such that
$$|\{k \in \{0,1,\ldots,p-1\}: f_p(z) \neq 0 \enskip {\text for \, all} \enskip z \in I_{k,\eta}\}| \geq (1-\varepsilon)p$$
for all sufficiently large primes $p \geq N_{\varepsilon}$.
\endproclaim

\head 4. Proofs of the Lemmas \endhead

\demo{Proof of Lemma 3.2} Let
$$Q(z) = c\prod_{j=1}^m{(z-a_j)}\,, \qquad c,a_j \in {\Bbb C}\,,$$
with some $m \leq p$. Without loss of generality we may assume that $c=1$. Note that
$$|a_j^p-1|^{1/p} \leq (2|a_j|^p)^{1/p} = 2^{1/p}|a_j|\,, \qquad |a_j| \geq 1\,,$$
while
$$|a_j^p-1|^{1/p} \leq 2^{1/p}\,, \qquad |a_j| < 1\,.$$
Multiplying these inequalities for $j=1,2,\ldots, m$, we obtain
$$\left(\prod_{j=0}^{p-1}{|Q(\zeta_p^j)|}\right)^{1/p} =
\left(\prod_{j=0}^m{|a_j^p-1|}\right)^{1/p} \leq
2^{m/p}\prod_{j=0}^m{\max\{|a_j|,1\}} \leq 2M_0(Q)\,.$$
\qed \enddemo

\demo{Proof of Lemma 3.3} Let
$$Q(z) = c\prod_{j=1}^m{(z-a_j)}\,, \qquad c,a_j \in {\Bbb C}\,,$$
with some $m \leq p$, where $a_j = b_j$, $j=1,2,\ldots,k$.
Without loss of generality we may assume that $c=1$. 
Note that
$$|a_j^p-1|^{1/p} \leq 2^{1/p}\max\{|a_j|,1\}\,, \qquad j=k+1,k+2,\ldots,m\,,$$
and
$$|a_j^p-1|^{1/p} \leq \left(2\cos \frac{\eta}{2}\right)^{1/p} = 
\left(2\cos \frac{\eta}{2}\right)^{1/p}|a_j|\,, \quad j=1,2,\ldots,k\,,$$
Multiplying these inequalities for $j=1,2,\ldots, m$, we obtain
$$\split \left(\prod_{j=0}^{p-1}{|Q(\zeta_p^j)|}\right)^{1/p} =
\left(\prod_{j=0}^m{|a_j^p-1|}\right)^{1/p} & \leq 2^{m/p}\left(\cos \frac{\eta}{2} \right)^{k/p} \prod_{j=0}^m{\max\{|a_j|,1\}} \cr 
& \leq 2 \left(\cos \frac{\eta}{2} \right)^{k/p} M_0(Q)\,. \cr \endsplit$$
\qed \enddemo

\demo{Proof of Lemma 3.7}
By Lemma 3.1 If $\zeta_p = e^{2i\pi/p}$ then, for all $k$ not divisible by $p$ we have  $|f_p(\zeta_p^k)| = \sqrt{p}$, and 
hence $|F_p(\zeta_p^k)| = \sqrt{p}$. Moreover 
$$F_p(\zeta_p^k) = (\zeta_p^k)^{-p/2} \sum_{a=1}^{p-1} {\left( \frac ap \right){\zeta_p^{ak}}} = 
(-1)^k \left( \frac kp \right) \sum_{a=1}^{p-1} {\left( \frac{ak}{p} \right){\zeta_p^{ak}}} = 
(-1)^k \left( \frac kp \right)F_p(\zeta_p)\,.$$
Therefore if $\displaystyle{\left( \frac kp \right) = \left( \frac {k+1}{p} \right)}$, then $\displaystyle{H_p\left(\frac{k}{p}\right)}$ and 
$\displaystyle{H_p\left(\frac{k+1}{p}\right)}$ have different signs. Since $H_p(t)$ is real-valued and continuous on the real line, 
it must have a zero $k/p$ and $(k+1)/p$ by the Intermediate Value Theorem. 
However, by Lemma 2 in [C-00] we have 
$$\left|\left\{k \in \{1,2,\ldots,p-2\}: \enskip \left( \frac kp \right) = \left( \frac {k+1}{p} \right)\right\}\right| = \frac{p-3}{2}\,,$$ 
and hence the values of $k \in \{0,1,\ldots,p-1\}$ for which $H_p$ has a zero between $k/p$ and $(k+1)/p$ is at least $(p-3)/2$.
\qed \enddemo

\demo{Proof of Lemma 3.9}
Note that 
$$I_\delta := \int_0^{\infty} {\sin(\delta \pi x) C(x) \, \frac{dx}{x}}$$
converges for every fixed $\delta > 0$, and  
$$\lim_{\delta \rightarrow 0+}{I_{\delta}} = 0\,.$$
Indeed, there is an absolute constant $c_1 > 0$ such that 
$$C(x) \leq c_12^{-3x/\pi}\,, \qquad x \geq 1\,,$$
as 
$$\left|\cos \left(\frac{2x}{2k+1}\right)\right| < \frac 12\,, \qquad \frac{3x}{\pi} < 2k+1 < \frac{6x}{\pi}\,.$$ 
Also, 
$$\left|\frac{\sin(\delta \pi x)}{x}\right| \leq \delta \pi\,, \qquad x > 0\,.$$
Therefore 
$$I_\delta \leq \int_0^{\infty} {\left| \frac{\sin(\delta \pi x)}{x} \right| \, |C(x)| \, dx} \leq A_{\delta} + B_{\delta}\,,$$
where
$$A_{\delta} := \int_0^{\delta^{-1/2}}{\left|\frac{\sin(\delta \pi x)}{x}\right| \, |C(x)|\,dx} 
\leq \delta^{-1/2}\delta \pi \leq \delta^{1/2}\pi \,,$$ 
and
$$B_{\delta} := \int_{\delta^{-1/2}}^{\infty}{\frac{|C(x)|}{x}\,dx} \leq \delta^{1/2} \int_{\delta^{-1/2}}^{\infty}{c_12^{-3x/\pi} \, dx} 
\leq \delta^{1/2} \frac{c_1\pi}{3\log 2}\,.$$
So by choosing $\delta > 0$ so that
$$I_{\delta} \leq A_{\delta} + B_{\delta} \leq \delta^{1/2}\pi + \delta^{1/2} \frac{c_1\pi}{3\log 2} \leq \frac{\pi \varepsilon}{2}\,,$$ 
the lemma follows from Lemma 3.8. 
\qed \enddemo

\demo{Proof of Lemma 3.10}
Suppose there are $0 \leq k_1 < k_2 < \cdots < k_m \leq p-1$ such that  
$$t_j \in I_{k_j}\,, \qquad |f_p^\prime(t_j)| \geq \gamma p^{3/2}\,, \qquad j=1,2,\ldots,m\,.$$  
Then
$$0 < t_1 < t_2 < \cdots < t_m < 2\pi\,,$$
and
$$\delta :=\min \left\{t_2 - t_1, t_3 - t_2,\dots, t_m - t_{m-1}, 2\pi -\left(t_m - t_1 \right) \right\} \geq \frac{\pi}{p}\,.$$
Hence by the large sieve inequality formulated in Lemma 3.6 and the Parseval formula appliied to 
$g_p(z) := z^{(3-p)/2}f_p^\prime(z)$ we get
$$\split m\gamma^2 p^3 & \leq \sum_{j=1}^{m}\left| f_p^\prime \left( e^{it_j} \right) \right|^{2} 
= \sum_{j=1}^{m}\left| g_p \left( e^{it_j} \right) \right|^{2} \cr 
& \leq \left(\frac{2(p-1)/2 + 1}{2\pi}+\delta ^{-1}\right) \int_{0}^{2\pi }{| g_p(e^{it})|^{2}\,dt} \cr 
& = \left(\frac{2(p-1)/2 + 1}{2\pi}+\delta ^{-1}\right) \int_{0}^{2\pi }{|f_p^\prime (e^{it})|^{2}\,dt} \cr 
& \leq \left(\frac{p}{2\pi} + \frac{p}{\pi} \right) 2\pi  \, \frac{(p-1)p(2p-1)}{6} = 3p \, \frac{(p-1)p(2p-1)}{6} \cr 
& \leq \frac{p^4}{2}\,. \cr \endsplit$$
\qed \enddemo

\demo{Proof of Lemma 3.11} 
Let $\varepsilon >0$. By Lemma 3.9 there is a $\delta > 0$ depending only on $\varepsilon >0$ such that 
$$\left|\left\{k \in \{1,2,\ldots,p\}: |f_p(e^{i(2k+1)\pi/p})| > \delta \sqrt{p}\right\}\right| \geq (1-\varepsilon/2)p \tag 4.1$$
for all sufficiently large primes $p \geq N_{\varepsilon}$.
Let $\gamma := \varepsilon^{-1/2}$. By Lemma 3.10 we have 
$$\left|\left\{k \in \{0,1,\ldots,p-1\}: \max_{z \in I_{k,\pi/2}} {|f_p^{\prime}(z)|} \leq \gamma p^{3/2}\right\}\right| \geq p - \gamma^{-2}p/2 
= (1-\varepsilon/2)p\,. \tag 4.2$$
Now let 
$$A_{p,\delta,\gamma} := \left\{k \in \{1,2,\ldots,p\}: |f_p(e^{i(2k+1)\pi/p})| > \delta \sqrt{p}\,, \quad 
\max_{z \in I_{k,\pi/2}} {|f_p^{\prime}(z)|} \leq \gamma p^{3/2}\right\}\,.$$
By (4.1) and (4.2) we obtain
$$|A_{p,\delta,\gamma}| \geq (1-\varepsilon)p\,. \tag 4.3$$
Let $0 < \eta < \min\{\delta/\gamma,\pi/2\}$. Observe that $k \in A_{p,\delta,\gamma}$ implies that 
$f_p$ does not vanish in  $I_{k,\eta}$. Indeed, $z := e^{it} \in I_{k,\eta}$ implies
$$\left|t - \frac{(2k+1)\pi}{p}\right| < \frac{\eta}{p}\,,$$
and hence 
$$\split |f_p(z)| \geq & \, |f_p(e^{i(2k+1)\pi/p})| - |f_p(z) - f_p(e^{i(2k+1)\pi/p})| \cr   
> & \, \delta \sqrt{p} - \int_{(2k+1)\pi/p}^t{f_p^\prime(e^{i\tau}) e^{i\tau} \, d\tau} \cr 
\geq & \, \delta \sqrt{p} - \int_{(2k+1)\pi/p}^t{|f_p^\prime(e^{i\tau})| |e^{i\tau}| \, d\tau} 
\geq \delta \sqrt{p} - \frac{\eta}{p}\,\gamma p^{3/2} \cr 
\geq & \,\delta \sqrt{p} - \delta \sqrt{p} = 0 \cr \endsplit$$ 
for all sufficiently large primes $p \geq N_{\varepsilon}$, and the lemma follows from (4.3).
\qed \enddemo

\head Proof of Theorem 2.1 \endhead

Now we are ready to prove the theorem.

\demo{Proof of Theorem 2.1} 
As in Lemma 3.11 let the subarcs $I_{k,\eta}$ of the unit circle ${\partial D}$ be defined by    
$$I_{k,\eta} := \left\{e^{it}: \left|t - \frac{(2k+1)\pi}{p}\right| < \frac{\eta}{p}\right\}\,, \qquad k=0,1,\ldots,p-1\,.$$
It follows from Lemma 3.11 that for $\varepsilon := 1/8$ there is an $\eta > 0$ such that
$$|\{k \in \{0,1,\ldots,p-1\}: f_p(z) \neq 0 \enskip \text {for all} \enskip z \in I_{k,\eta}\}| \geq i\frac{7p}{8}$$
for all sufficiently large primes $p$. Combining this with Lemma 3.7 we have that 
$$|\{k \in \{0,1,\ldots,p-1\}: f_p \enskip \text{has a zero on} \enskip I_{k,\pi} \enskip \text{and} \enskip 
f_p(z)\neq 0 \enskip \text{for all} \enskip z \in I_{k,\eta}\}| \geq p/4$$ 
for all sufficiently large primes $p$. Hence the assumptions of Lemma 3.3 are satisfied with 
$Q := f_p$ and $k \geq p/4$ for all sufficiently large primes.
Suppose that $1$ is a zero of $f_p$ with multiplicity $m = m(p)$.
By either Lemma 3.4 or Lemma 3.5 we have $m = O(p^{1/2})$.
Let $g_p(z) := f_p(z)/h_m(z)$ with $h_m(z) := (z-1)^m$.
Note that $|g_p(1)|$ is a nonzero integer, hence $|g_p(1)| \geq 1$.
Also, $h_m$ is monic and has all its zeros on the unit circle, hence
$M_0(h_m)=1$. Combining these with
the multiplicative property of the Mahler measure, Lemma 3.3 applied to $g_p$ with $k \geq p/4$, Lemma 3.1, and the fact that 
$m = O(p^{1/2})$ implies that
$$\lim_{p \rightarrow \infty}{p^{(-(1/2)+m)/p}} = 1\,,$$
we conclude that there are absolute constants 
$$c_1 := \left( 2\left(\cos \frac{\eta}{2} \right)^{1/4} \right)^{-1} > c_2 > \frac 12$$
such that   
$$\split M_0(f_p) & = M_0(g_p)M_0(h_m) = M_0(g_p) \cr
& \geq \left( 2\left(\cos \frac{\eta}{2} \right)^{k/p} \right)^{-1} \left(|g_p(1)|\prod_{k=1}^{p-1}{|g_p(\zeta_p^k)|}\right)^{1/p} \cr 
& \geq c_1 \left(|g_p(1)|\prod_{k=1}^{p-1}{\left|\frac{f_p(\zeta_p^k)}{(\zeta_p^k-1)^m}\right|}\right)^{1/p}  
\geq c_1 \left(\prod_{k=1}^{p-1}{\left|\frac{f_p(\zeta_p^k)}{(\zeta_p^k-1)^m}\right|}\right)^{1/p} \cr 
& = c_1 \frac{(p^{1/2})^{(p-1)/p}}{p^{m/p}} = c_1 p^{(p-1)/(2p) - m/p} = c_1 p^{1/2} p^{-((1/2)+m)/p} \cr 
& \geq c_2\sqrt{p} \cr \endsplit$$
for all sufficiently large primes $p$. 
\qed \enddemo

\subhead {5. Acknowledgment}\endsubhead
The author thanks Stephen Choi for his careful reading of the paper 
and for his comments.

\enddocument